\documentclass[a4paper,12pt]{amsart}

\def\d{\partial}
\def\CP1{\mathbb{C}\mathrm{P}^1}
\def\C{\mathbb{C}}
\def\ad{\mathrm{ad}\,}

\newtheorem{theorem}{Theorem}
\newtheorem{proposition}[theorem]{Proposition}

\title[Goulden-Jackson-Vakil formula]{On the structure of Goulden-Jackson-Vakil formula}

\author{S. Shadrin}

\address{Korteweg-de~Vries Institute for Mathematics, University of Amsterdam,
Plantage Muidergracht 24, 1018 TV Amsterdam, The Netherlands}

\email{s.shadrin@uva.nl}

\address{Department of Mathematics, Institute of System Research, Nakhimovsky prospekt 36-1, Moscow 117218, Russia}

\email{shadrin@mccme.ru}

\begin{document}

\begin{abstract}
We study the structure of the Goulden-Jackson-Vakil formula that relates Hurwitz numbers to 
some conjectural ``intersection numbers'' on a conjectural family of varieties $X_{g,n}$ of dimension $4g-3+n$. We give explicit formulas for the properly arranged generating function for these ``intersection numbers'', and prove that it satisfies Hirota equations. This generalizes and substantially simplifies our earlier results with Zvonkine.
\end{abstract}

\maketitle

\tableofcontents

\section{Introduction}

\subsection{Hurwitz numbers}

Consider a ramified covering of degree $d$ over the sphere $\CP1$ by a smooth Riemann surface of genus $g$. 
Assume that there is a total ramification over $0\in\CP1$ (that is, the monodromy at $0$ is a cycle of length $d$);
there are $n$ preimages of $\infty\in\CP1$ with the multiplicities $b_1,\dots,b_n$; and there is exactly one simple ramification over each of $m=2g-1+n$ fixed points $z_1,\dots,z_m\in\CP1\setminus\{0,\infty\}$ (the number $m$ is determined by the Riemann-Hurwitz formula). 

There is a finite number of such coverings. The Hurwitz number $h_{g,b_1,\dots,b_n}$ counts these coverings with weights; the weight of a covering equals to the inversed order of its automorphism group.

\subsection{GJV-formula}

Goulden, Jackson, and Vakil studied the structure of these Hurwitz numbers in~\cite{GoJaVa05}. They observed that these numbers can be represented as
\begin{equation}
\frac{h_{g,b_1,\dots,b_n}}{d\cdot m!}
=
\int_{X_{g,n}}
\frac{1-\lambda_2+\cdots\pm\lambda_{2g}}{(1-b_1\psi_1)\cdots(1-b_n\psi_n)}.
\end{equation}

Here $X_{g,n}$ is a conjectural complex algebraic variety of complex dimension $4g-3+n$, and $\psi_1,\dots,\psi_n\in H^2(X_{g,n},\C)$ and $\lambda_{2i}\in H^{4i}(X_{g,n},\C)$, $i=1,\dots,g$, are some conjectural cohomology classes on it. There is a reasonable hope that $X_{g,n}$ might be something like a suitable compactification of the universal Picard variety; and $\psi$- and $\lambda$-classes might play the same significant role as the usual $\psi$- and $\lambda$- classes do in geometry of the moduli space of curves. 

\subsection{Relation to Hirota equations}

A consequence of the GJV conjecture is that the ``intersection numbers'' 
\begin{equation}
\langle \lambda_{2i}\tau_{d_1}\cdots\tau_{d_n}\rangle:=\int_{X_{g,n}}\lambda_{2i}\psi_1^{d_1}\cdots\psi_n^{d_n},
\end{equation}
$4g-3+n=2i+d_1+\cdots d_n$, are not only some important combinatorial constans (defined via Hurwitz numbers and GJV-formula), but they might also be related to geometry and deserve further study.

The first step in this direction was done by me and Zvonkine in~\cite{ShZv07}. We proved that the generating series for the numbers $\langle \tau_{d_1}\cdots\tau_{d_n}\rangle$ is a solution to Hirota hierarchy and to the linearized KP hierarchy simulteneously. This is a complete analogue of the Witten-Kontsevich theorem~\cite{ChLiLi06,Ka07,KaLa06,KiLi05,Ko02,Mi07,OkPa01,VeVe91,Wi91} for the intersection numbers of $\psi$-classes on the moduli space of curves. 

\subsection{Kazarian's study of Hodge integrals}

In the case of the usual moduli space, there is another formula that relates the intersection number with $\psi$-classes and at most one $\lambda$-class (Hodge integrals) to the combinatorics of ramified coverings. It is the celebrated ELSV-formula~\cite{EkLaShVa01,GrVa03}. It was recently used by Kazarian~\cite{Ka07} in order to give a generalization of the Witten-Kontsevich theorem. He rearranges the generating series for Hodge integrals in such a way that is appears to be a solution of Hirota hierarchy. 

In this note we apply the technique proposed by Kazarian in order to expand the structure of GJV-formula. It works perfectly; we also relate the generating series for our ``intersection numbers'' $\langle\lambda_{2i}\tau_{d_1}\cdots\tau_{d_n}\rangle$ to Hirota equations. Moreover, all computations in this case appear to be much more simple. In particular, we manage to give a very explicit formulas for the generating series of our conjectural ``intersection numbers''. Of course, we hope that this computation will be helpful for the search of the proper family of varieties $X_{g,n}$ of dimension $4g-3+n$, whose existence and explicit description is demanded as a part of the Goulden-Jackson-Vakil conjecture.

\subsection{Hirota equations}

In this paper we call by Hirota equation the bilinear Hirota form of the Kadomtsev-Petviashvili hierarchy. So, a formal power series satisfies Hirota equations if and only if it is a tau-function of the KP hierarchy.
For the complete definition and discussion of Hirota equations, we refer the reader to~\cite{DaKaJiMi83, KaRa87, Ka07, MiJiDa00, ShZv07}. For our purposes in this paper, we use a small number of properties of Hirota equations. That is:
\begin{enumerate}
\item For any $c\in\C$, the function $c+q_1$ satisfies Hirota equations.
\item The operators 
\begin{align*}
\Lambda_1& =\sum_{i=2}^\infty q_i\frac{(i-1)\d}{\d q_{i-1}}, \\
\Lambda_a &=\sum_{i=1}^\infty q_i\frac{(i-a)\d}{\d q_{i-a}}, \qquad a\leq 0 \\
M_0&=\frac{1}{2}\sum_{i,j=1}^\infty q_iq_j\frac{(i+j)\d}{\d q_{i+j}}+
\frac{1}{2}\sum_{i,j=1}^\infty q_{i+j}\frac{ij\d^2}{\d q_i\d q_j} \\
M_1&=\frac{1}{2}\sum_{i+j\geq 1}^\infty q_iq_j\frac{(i+j-1)\d}{\d q_{i+j-1}}+
\frac{1}{2}\sum_{i,j=1}^\infty q_{i+j+1}\frac{ij\d^2}{\d q_i\d q_j} \\
M_2&=\frac{1}{2}\sum_{i+j\geq 2}^\infty q_iq_j\frac{(i+j-2)\d}{\d q_{i+j-2}}+
\frac{1}{2}\sum_{i,j=1}^\infty q_{i+j+2}\frac{ij\d^2}{\d q_i\d q_j} 
\end{align*}
are the infinitesimal symmetries of Hirota equations.
\item Hirota equations are preserved by the rescaling of variables $q_i\leftrightarrow u^i q_i$ (simulteneously for all $i=1,2,\dots$; $u$ is a formal parameter).
\end{enumerate}

Let us give an example of the argument that we use below. Since $c+q_1$ satisfies Hirota equations,
the series $\exp(\Lambda_1)(c+q_1)=c+\sum_{i=1}^\infty q_i$ also satisfies Hirota equations.

An important point for us is that if both $f$ and $1+f$ satisfy Hirota equations, it means that $f$ also satisfies linearized KP hierarchy. So we keep track of an opportunity to add an arbitrary constant to a solution of Hirota equations through the text.

\subsection{Acknoledgements}

I am grateful to M.~Kazarian for the helpful remarks.

\section{Rearranged generating series and Hirota equations}

In this section, we study the generating series for all integrals involved in Goulden-Jackson-Vakil formula.

\subsection{Hirota equations and an explicit formula}
We define a sequence of linear functions $T_k$, $k=0,1,\dots$, in formal variables $q_i$, $i=1,2,\dots$. We set $T_0=q_1$, and $T_{k+1}=(u\Lambda_0+\Lambda_1) T_k$. We list the first few expressions:
\begin{align}\label{eq:T}
T_0 &= q_1,\\
T_1 &= uq_1 + q_2, \notag \\
T_2 &= u^2 q_1 + 3uq_2 + 2q_3, \notag\\
T_3 &= u^3 q_1 + 7u^2q_2 + 12uq_3 + 6q_4,\notag
\end{align}
and so on. 

Consider the generating series for the intersection numbers with $\lambda$-classes defined as
\begin{equation}
G(u,q_1,q_2,\dots)=\sum_{j,k_1,k_2,\dots} (-1)^{j}\langle \lambda_{2j} \tau_0^{k_0} \tau_1^{k_1} \dots\rangle u^{2j+1}
\frac{T_0^{k_0}}{k_0!} \frac{T_1^{k_1}}{k_1!} \dots.
\end{equation}
In particular, we see that that the expansion of $G$ in $u$ starts with the terms 
\begin{equation}
G=uF+\frac{u^2}{2}\Lambda_{-1} F+\dots,
\end{equation}
where $F$ is the generation function for the intersection numbers without $\lambda$-classes:
\begin{equation}
F(t_0,t_1,\dots)=\sum_{k_0,k_1,\dots} \langle \tau_0^{k_0} \tau_1^{k_1} \dots\rangle 
\frac{(0!q_1)^{k_0}}{k_0!} \frac{(1!q_2)^{k_1}}{k_1!} \dots.
\end{equation}

\begin{theorem}\label{thm1} For any function $c=c(u)$, the series 
\begin{equation}\label{eq:fullseries}
\tau = c(u)+\frac{q_1+q_1q_2}{u}+q_1^2+\left(\Lambda_0+\frac{1}{u}\Lambda_1\right)^2 G(u,q_1,q_2,\dots)
\end{equation}
is a solution of the Hirota equations in variables $q_i$, $i=1,2,\dots$ ($u$ is just a parameter).
\end{theorem}

In fact, we can give an explicit formula for $\tau$.
\begin{theorem}\label{thm2} We have:
\begin{equation}\label{eq:tau2}
\tau=\exp(M_2+2uM_1+u^2M_0)\left(c(u)+\frac{q_1}{u}\right).
\end{equation}
\end{theorem}

\subsection{Generating series for Hurwitz numbers}

Consider the generating series for Hurwitz numbers,
\begin{align}
H(\beta,p_1,p_2,\dots)& :=\sum_{g,n} H_{g,n} \\
& :=\sum_{g,n}  \frac{1}{n!}\sum_{b_1,\dots,b_n}\frac{h_{g,b_1,\dots,b_n}}{b\cdot m!}p_{b_1}\cdots p_{b_n}\beta^m.
\notag
\end{align}
It is proved in~\cite{ShZv07} (following the observations made before in~\cite{Ok00,KaLa06}) that $c(\beta)+\Lambda_0^{2}H$ is a solution of the Hirota equations in variables $p_i$, $i=1,2,\dots$, for any function $c(\beta)$.

For completeness, let us remind the proof. Hurwitz numbers satisfy a so-called cut-and-join equation (see, e.~g.~\cite{GoJaVa05}). This practically means that $H=\exp(\beta M_0)H_{0,1}$. Since $[M_0,\Lambda_0]=0$, we have:
\begin{equation}
\Lambda_0^2H=\exp(\beta M_0)\Lambda_0^2H_{0,1}=\exp(\beta M_0)\sum_{i=1}^\infty p_i.
\end{equation}
Since the series $c(\beta)+\sum_{i=1}^\infty p_i$ satisfies the Hirota equations, and $M_0$ is an infinitesimal symmetry of the Hirota equations, we conclude that $c(\beta)+\Lambda_0^2H$ also satisfies the Hirota equations.

\subsection{Rearranging of GJV-formula}

The GJV-formula is applied to all components of $H$ except for $H_{0,1}$ and $H_{0,2}$. Denote by $H_n$ the sum $n!\sum_{g=0}^n H_{g,n}$ for $n\geq 3$ and the sum $n!\sum_{g=1}^n H_{g,n}$ for $n=1,2$. That is, 
\begin{equation}
H=\sum_{n=1}^\infty \frac{H_n}{n!}+H_{0,1}+H_{0,2}.
\end{equation}

Using that $m=\dim X_{g,n}/2+n/2+1/2$, we obtain:
\begin{align}
H_n &=\sum_{g,b_1,\dots,b_n}\int_{X_{g,n}}\frac{1-\lambda_2+\dots\pm\lambda_{2g}}{(1-b_1\psi_1)\cdots (1-b_n\psi_n)}p_{b_1}\cdots p_{b_n}\beta^{m} \\
&= \beta\sum_{g,b_1,\dots,b_n}\int_{X_{g,n}}(1-\beta\lambda_2+\dots\pm\beta^g\lambda_{2g})\prod_{i=1}^n
\frac{\beta^{1/2}p_{b_i}}{(1-\beta^{1/2}b_i\psi_i)} \notag \\
& = \beta\left\langle(1-\beta\lambda_2+\beta^2\lambda_4-\dots)\prod_{i=1}^n (\sum_{d\geq 0} \tau_d T_d)\right\rangle, \notag
\end{align}
where
\begin{equation}
T_d=\sum_{b\geq 1} \beta^{1/2}p_b\cdot b^d\beta^{d/2}=\beta^{(d+1)/2}\sum_{b\geq 1} b^d p_b.
\end{equation}
There is an obvious way to define $T_d$ recursively. We just set $T_0=\beta^{1/2}\sum_{b\geq 1}p_b$, and $T_{k+1}=\beta^{1/2}\Lambda_0 T_k$.

\subsection{Change of variables}

We need the following change of variables. First, we replace $\beta^{1/2}$ by $u$. Then we rescale the variables by setting $p_b=q_b/u^b$. And then we replace $q_i$ with $\exp(-\Lambda_1/u)q_i$. 

In other words, it is a linear triangular change of variables given by 
\begin{equation}\label{eq:change}
p_b=\sum_{i=b}^\infty \frac{1}{u^i}(-1)^{i-b}\binom{i-1}{b-1}q_i.
\end{equation}
The same change of variables (in a bit different notation) is used in~\cite{ShZv07}.
Under this change of variable a series $f(\beta^{1/2},p_1,p_2,\dots)$ transforms into $g(u,q_1,q_2,\dots):=\exp(-\Lambda_1/u)f(u,q_1/u,q_2/u^2,\dots)$.

Another way to do the same is to replace $p_b$ with $\exp(-\Lambda_1)p_b$ and then to rescale $p_b=q_b/u^b$.
Note that $\Lambda_1$ in an infinitesimal symmetry of Hirota equations. The rescaling $p_b=q_b/u^b$ also preserves them. Therefore, this change of variables preserves the property to be a solution to Hirota equations. 

\subsection{Proof of Theorems~\ref{thm1} and~\ref{thm2} }

We have already observed that $c(\beta)+\Lambda_0^2H=\exp(\beta M_0)(c(\beta)+\sum_{i=1}^\infty p_i)$ is a solution to Hirota equations.
Under the change of variables~\eqref{eq:change}, $c(\beta)+\sum_{i=1}^\infty p_i$ turns into $c(u)+q_1/u$ (it is an explicit computation). Since $[M_0,\Lambda_1]=2M_1$, $[M_1,\Lambda_1]=M_2$, and $[M_2,\Lambda_1]=0$, the operator $u^2M_0$ turns into 
\begin{equation}
\exp(-\Lambda_1/u)u^2M_0\exp(\Lambda_1/u)=u^2M_0+2uM_1+M_2.
\end{equation}
This implies that under the change of variables~\eqref{eq:change}, $c(\beta)+\Lambda_0^2H$ turns into
\begin{equation}\label{eq:tau}
\tau=\exp(M_2+2uM_1+u^2M_0)\left(c(u)+\frac{q_1}{u}\right).
\end{equation}

On the other hand, under the change of variables~\eqref{eq:change}, $\Lambda_0^2H_{0,1}$ turns into $q_1/u$, $\Lambda_0^2H_{0,2}$ turns into $q_1q_2/u+q_1^2$, and $T_0$ turns into $q_1$ (all these observations are simple explicit computations). Also, since $[\Lambda_0,\Lambda_1]=\Lambda_1$, it follows that $\beta^{1/2}\Lambda_0$ turns into $\exp(-\Lambda_1/u)u\Lambda_0\exp(\Lambda_1/u)=u\Lambda_0+\Lambda_1$. 

Therefore, $H_n$ turns into $\sum_{j,b_1,\dots,b_n}(-1)^{j}\langle\lambda_{2j}\tau_{b_1}\cdots\tau_{b_n}\rangle u^{2j+1}T_{b_1}\cdots T_{b_n}$, where $T_b$ is defined in variables $q$ as in equation~\eqref{eq:T}. 

So, $c(\beta)+\Lambda_0^2H=c(\beta)+\Lambda_0^2H_{0,1}+\Lambda_0^2H_{0,2}+\Lambda_0^2\left(\sum_{n\geq 1} H_n/n!\right)$ turns into 
$c(u)+q_1/u+q_1q_2/u+q_1^2+(\Lambda_0+\Lambda_1/u)^2 G$. The last expression must coincide with $\tau$ in equation~\eqref{eq:tau}, and this completes the proof of both theorems.


\section{Intersections of $\psi$-classes}

In this section, we study the formal power series 
\begin{equation}
F(q_1,q_2,\dots)=\sum_{k_0,k_1,\dots} \langle \tau_0^{k_0} \tau_1^{k_1} \dots\rangle 
\frac{(0!q_1)^{k_0}}{k_0!} \frac{(1!q_2)^{k_1}}{k_1!} \dots.
\end{equation}

\subsection{An explicit formula}

Using two different expressions for $\tau$ given by formulas~\eqref{eq:fullseries} and~\eqref{eq:tau2}, we see that
\begin{equation}\label{eq:lambda1}
\Lambda_1^2F+q_1+q_1q_2=\exp(M_2)q_1.
\end{equation}

It was proved in~\cite{GoJaVa05} (and we reprove this below) that $F$ satisfies the string equation, that is,
\begin{equation}\label{eq:string}
\frac{\d}{\d q_1} F = \Lambda_1 F +\frac{q_1^2}{2}.
\end{equation}
It implies that 
\begin{equation}
\Lambda_1^2 F+ q_1+q_1q_2 = \frac{\d^2 F}{\d q_1^2}.
\end{equation}

Thus we prove the following theorem.
\begin{theorem} We have:
\begin{equation}
\frac{\d^2}{\d q_1^2}F=\exp(M_2)q_1
\end{equation}
\end{theorem}

In particular, an obvious corollary of this explicit formula is Theorem 2 in~\cite{ShZv07}:
\begin{theorem} For any $c\in\C$, the series $c+\d^2 F/\d q_1^2$ satisfies Hirota equations.
\end{theorem}

\subsection{String equation}

Let us prove the string equation~\eqref{eq:string}. Since $[\d/\d q_1,M_2]=\Lambda_1$, and $\Lambda_1$ commutes with $M_2$, we have:
\begin{equation}
\frac{\d \exp(M_2)q_1}{\d q_1}=\Lambda_1\exp(M_2)q_1+1.
\end{equation}
If we substitue in this formula $\exp(M_2)q_1$ with $\Lambda_1^2F+q_1+\Lambda_1(q_1^2/2)$ (using equation~\eqref{eq:lambda1}), and use that $\d/\d q_1$ commutes with $\Lambda_1$, we obtain
\begin{equation}
\Lambda_1^2\left(\frac{\d F}{\d q_1}\right)=\Lambda_1^2 \left(\Lambda_1 F+ \frac{q_1^2}{2}\right).
\end{equation}
This implies the string equation~\eqref{eq:string}.

\subsection{Towards Virasoro constrains}

Though we already have a nice closed formula for $\d^2 F/\d q_1^2$, it is still might be interesting to find an analog of Virasoro contrains for it. Indeed, our goal is to find a family of varieties with the intersection theory controlled by $F$. And we know that in the usual case of the moduli space of curves the geometry of degenerations of curves is intimately related to the existence and the particular form of Virasoro constrains, see~\cite{Mi07}. However, a part of the problem is that $F$ is not an exponential generating function; and in our case it doen't make any sense to exponentiate it.

We prove a sequence of some strange equations for $\exp(M_2)q_1$; and the first two of them are indeed the string equation and the dilaton equation, already proved in~\cite{GoJaVa05}.

\begin{proposition} For any $n\geq 1$,
\begin{equation}
\frac{n \d }{\d q_n}\exp(M_2)q_1=\Lambda_{2-n}\exp(M_2)q_1+\exp(M_2)q_1^{n-1}.
\end{equation}
\end{proposition}

\subsubsection{Proof}
The proof is based on several observations. Denote by $\ad y(x)$ the operator $[x,y]$. Then the operator $n\d/\d q_n$ is equal to the sum 
\begin{equation}
\frac{n\d}{\d q_n}=\sum_{i=0}^\infty \exp(-\ad M_2) \frac{\ad M_2^i}{i!}\left(\frac{n\d}{\d q_n}\right).
\end{equation}

We show below that among all operators
\begin{equation}
O_i:=\exp(-\ad M_2) \frac{(\ad M_2)^i}{i!}\left(\frac{n\d}{\d q_n}\right)
\end{equation}
only $O_{n-1}$ and $O_n$ act nontrivially on $\exp(M_2)q_1$. Moreover, one can prove that $O_{n-1}\exp(M_2)q_1=n\exp(M_2)q_1^{n-1}$.

Observe that $n\Lambda_{2-n}=[\ad M_2, n\d/\d q_n]$ is also equal to $\sum_{i=1}^\infty iO_i$. Thus we have: 
\begin{align}
\frac{n\d}{\d q_n} \exp(M_2)q_1 & = n\exp(M_2)q_1^{n-1}+O_n\exp(M_2)q_1 \\
\Lambda_{2-n}\exp(M_2)q_1 & =(n-1)\exp(M_2)q_1^{n-1}+O_n\exp(M_2)q_1, \notag
\end{align}
and this implies the statement of the theorem.

\subsubsection{Action of $O_i$}
It is a simple observation that 
\begin{equation}
i!O_i\exp(M_2)q_1=\exp(M_2)\left(\left(\ad M_2\right)^i\left(\frac{n\d}{\d q_n}\right)\right) q_1.
\end{equation}
Then,
\begin{equation}
\left(\left(\ad M_2\right)^i\left(\frac{n\d}{\d q_n}\right) \right) q_1
=\left( \left(\ad M_2^{lin}\right)^i\left(\frac{n\d}{\d q_n}\right)\right) q_1,
\end{equation}
 where by $M_2^{lin}$ we denote the linear part of $M_2$. And then it is a straightforward calculation to 
show that for $i\geq 1$
\begin{multline}
\left(\ad M_2^{lin}\right)^i\left(\frac{n\d}{\d q_n}\right) = \\ \prod_{j=0}^{i-1}(n-j)\sum_{k_1,\dots,k_i} q_{k_1}\cdots q_{k_i}\frac{(k_1+\cdots+k_i+n-2i) \d}{\d q_{k_1+\cdots+k_i+n-2i}}
\end{multline}

Since $k_l\geq 1$, $l=1,\dots,i$, the terms with $\d/\d q_1$ appear only for $i\geq n-1$. But the coefficient $\prod_{j=0}^{i-1}(n-j)$ vanishes for $i\geq n+1$. So this operator can be applied non-trivially to $q_1$ only if $i=n-1$ or $n$. Then it is obvious that for $i=n-1$ it turns $q_1$ into $n!q_1^{n-1}$.

\end{document}